\documentclass[12pt]{amsart}
\usepackage[utf8]{inputenc}
\usepackage{fullpage}
\usepackage{graphicx}
\usepackage{float}
\usepackage{parskip}
\usepackage{amssymb,amscd,amsbsy,amsmath, amsthm, amsfonts, amsrefs}
\usepackage[OT2,OT1]{fontenc}
\usepackage{color}   %May be necessary if you want to color links
\usepackage[colorlinks=true,linkcolor=blue,citecolor=blue,urlcolor=blue]{hyperref}
\usepackage{algorithmicx} 
\usepackage{algpseudocode} 
\usepackage{comment}
\includecomment{A}
\newtheorem{theorem}{Theorem}%[section]

\newtheorem{lemma}[theorem]{Lemma}

\theoremstyle{definition}

\theoremstyle{remark}

\newfont{\cmbsy}{cmbsy10}
\newfont{\cmmib}{cmmib10}
\newcommand{\Orden}{\mathop{\hbox{\cmbsy O}}\nolimits}

\def\Z{\mathbf{Z}}

\def\R{\mathbf{R}}

\def\Hermite{\operatorname{\mathcal H}}

\begin{document}

\title{Explicit Formula and quasicrystal definition.}
\author{J. Arias de Reyna}
\address{Univ.~de Sevilla \\
Facultad de Matem\'aticas \\
c/Tarfia, sn
 \\
41012-Sevilla \\
Spain} 
%\thanks{Supported by  MINECO grant MTM2015-63699-P}
\email{arias@us.es}

\date{\today, \texttt{148-ZetaExplicit-v5.tex}}

\begin{abstract}

We show that the Riemann hypothesis is true if and only if the measure
\[\mu=-\sum_{n=1}^\infty\frac{\Lambda(n)}{\sqrt{n}}(\delta_{\log n}+\delta_{-\log n})+2\cosh(x/2)\,dx\]
is a tempered distribution. In this case it is the  Fourier transform of another measure 
\[\mathcal{F}\Bigl(\sum_{\gamma}\delta_{\gamma/2\pi}-2\vartheta'(2\pi t)\,dt\Bigr)=\mu.\]
We propose a definition of Fourier quasi-crystal to make sense of Dyson suggestion.
\end{abstract}

\maketitle

\section{Introduction}
Freeman Dyson in \cite{Dy} define a quasi-crystals as a distribution of discrete point masses whose Fourier transform is a distribution of discrete frequencies.
And then asserts: \emph{if the Riemann hypothesis is true, then the zeros of the zeta-function form a one dimensional quasi-crystal}. 
Recently Pavel Kurasov and Peter Sarnak \cite{KS} and  Alexander Olevskii and Alexander Ulanovskii \cite{OU} have determined all Fourier quasi-crystals with unit masses, i.e. all tempered measures 
\[\mu=\sum_{\lambda\in\Lambda}\delta_\lambda \quad
\text{with Fourier transforms of type}\quad
\widehat\mu=\sum_{s\in S}w_s\delta_s,\]
with $\Lambda$ and $S$ locally finite sets. 

Among these Fourier quasi-crystals there are none corresponding to the zeros of the zeta function.
Dyson  appear to refer to Weil's explicit formula, but one finds that Weil's explicit formula do not say exactly what Dyson claim. For example in  Delsarte  version \eqref{E:Delsarte3} we find  the term 
$\widehat f(i/2)$, so that even assuming the Riemann hypothesis one of the measures appear to have a support off the line.  The version of  Brad Rodgers in his thesis \cite{R}
\[\lim_{L\to\infty}\sum_{|\gamma|<L}\widehat{g}\Bigl(\frac{\gamma}{2\pi}\Bigr)-\int_{-L}^L\frac{\Omega(\xi)}{2\pi}\widehat{g}\Bigl(\frac{\xi}{2\pi}\Bigr)\,d\xi=
\int_{-\infty}^\infty[g(x)+g(-x)]e^{-x/2}\,d(e^x-\psi(e^x)),\]
is very near what Dyson have in mind but it is not worded as in Dyson and is true independent of the Riemann hypothesis.  Guinand in 1959 gave almost exactly what Dyson have in mind
If $g(x)=\sqrt{2/\pi}\int_0^\infty f(t)\sin(xt)\,dt$ and we assume the RH, then 
\begin{multline*}
\lim_{T\to\infty}\Bigl\{\sum_{\log n<T}\frac{\Lambda(n)}{\sqrt{n}}f(\log n)-\int_0^Te^{t/2}f(t)\,dt\Bigr\}-\frac12\int_0^\infty f(t)\Bigl(\frac{1}{t}-\frac{e^{-3t/2}}{\sinh t}\Bigr)\,dt\\=-\sqrt{2\pi}\lim_{T\to\infty}\Bigl\{\sum_{0<\gamma<T}g(\gamma)-\frac{1}{2\pi}\int_0^Tg(t)\log\frac{t}{2\pi}\,dt\Bigr\}.\end{multline*}

The object of this note is to give explicitly what I think Dyson had in mind. This is 
a trivial observation, but after so many years without anybody writing this explicitly I think it is not entirely superfluous.

In particular we prove that if the Riemann hypothesis is true there are two measures that are tempered distributions one being the Fourier transform of the other
\begin{equation}
\mathcal{F}\Bigl(\sum_\gamma\delta_{\gamma/2\pi}-2\vartheta'(2\pi t)\,dt\Bigr)=
-\sum_{n=1}^\infty\frac{\Lambda(n)}{\sqrt{n}}(\delta_{\log n}+\delta_{-\log n})+2\cosh(x/2)\,dx.
\end{equation}
We show also that the right hand side measure is a tempered distribution if and only if the Riemann hypothesis is true.

We may understand also the extra terms in our formulation. The measure
\[\sum_{n=1}^\infty\frac{\Lambda(n)}{\sqrt{n}}(\delta_{\log n}+\delta_{-\log n})\]
is not a tempered distribution because its mass in the interval $[-x,x]$ is approximately (see Lemma \ref{L:aux})
\[2\sum_{\log n\le x}\frac{\Lambda(n)}{\sqrt{n}}\sim 4e^{x/2}\sim\int_{-x}^x 2\cosh(y/2)\,dy.\]
Hence it is approximated by a measure of density $2\cosh(x/2)$.  This is what made possible that the difference may be a tempered distribution

Assuming the Riemann hypothesis the $\gamma_n$ are real and $\sum_\gamma\delta_{\gamma/2\pi}$ is a tempered distribution. The mass in $[-x,x]$ being \[2N(2\pi x)\sim 2\frac{\vartheta(2\pi x)}{\pi}=\int_{-x}^x2\vartheta'(2\pi t)\,dt \]
and it is natural to consider the measure $\sum_\gamma\delta_{\gamma/2\pi}-2\vartheta'(2\pi x)\,dx$. 

To make the words of Dyson true we may define a quasi-crystal as a discrete set of points $A$ in $\R$ such that the measure $\sum_{a\in A}\delta_{a}-\rho_A(t)\,dt$, where $\rho_A(t)$ is a continuous function representing the  mean density of points in  $A$ near $t$, is a tempered distribution with Fourier transform 
\[\mathcal{F}\Bigl(\sum_{a\in A} \delta_{a}-\rho_A(t)\,dt\Bigr)=\sum_{b\in B} w(b)\delta_{b} -\rho_B(x)\,dx.\]
where  $w(b)$ are some  weights and $\rho_B(x)$ is the density of the points $b\in B$ counted with his weights $w(b)$.

\section{Main theorem}

We will need the following version of the explicit formula 
\begin{theorem}[Delsarte]\label{T:Dels}
Let $F(s)$ be a holomorphic function on $-1\le \sigma\le 2$. Assume that for any $k\in N$,  $\lim_{t\to\infty}|t|^k \widehat{G}(\sigma+it)=0$ uniformly for $\sigma\in[-1,2]$. For $t$ in $\R$, put $f(t)=F(\frac12+it)$. Then 
\begin{equation}\label{E:Delsarte}
\sum_\rho F(\rho)-F(0)-F(1)=-\frac{1}{2\pi}\sum_{n=1}^\infty \frac{\Lambda(n)}{\sqrt{n}}\bigl(\widehat f(\tfrac{\log n}{2\pi})+\widehat f (-\tfrac{\log n}{2\pi})\bigr)-\frac{1}{2\pi}\int_{-\infty}^\infty f(t)\Phi(\tfrac12+it)\,dt,
\end{equation}
where $\rho$ run through the non-trivial zeros of $\zeta(s)$ repeated according to its multiplicity and 
\[\Phi(s)=\frac{\zeta'(s)}{\zeta(s)}+\frac{\zeta'(1-s)}{\zeta(1-s)}.\]
\end{theorem}

\begin{proof}
It is found in Delsarte \cite{D}. Its proof consist in applying Cauchy's Theorem to the function $F(s)\frac{\zeta'(s)}{\zeta(s)}$ in the strip $-1\le\sigma\le2$ using the functional equation and the Dirichlet series for $\zeta'(s)/\zeta(s)$.
\end{proof}

\begin{lemma}\label{L:aux} Assuming the Riemann hypothesis we have
\begin{equation}\label{E:201108-1}
\sum_{n\le x}\frac{\Lambda(n)}{\sqrt{n}}=2\sqrt{x}+\Orden(\log^3x).
\end{equation}
\end{lemma}
\begin{proof}
Riemann's hypothesis implies  (see Ingham \cite{I}*{Th.~30})
\[\psi(x)=\sum_{n\le x}\Lambda(n)=x+\Orden(x^{1/2}\log^2x).\]
Therefore by partial summation
\begin{align*}
\sum_{n\le x}\frac{\Lambda(n)}{\sqrt{n}}&=\frac{\psi(x)}{\sqrt{x}}+\frac12\int_1^x \psi(t)t^{-3/2}\,dt\\
&=\sqrt{x}+\Orden(\log^2x)+\frac12\int_1^x\frac{dt}{\sqrt{t}}+\frac12\int_1^x\frac{R(t)dt}{t^{3/2}}\\
&=2\sqrt{x}-1+\Orden(\log^2x)+\Orden(\log^3x).\qedhere
\end{align*}
\end{proof}

\begin{theorem}
The Riemann hypothesis is equivalent to say that the measure
\[\mu=-\sum_{n=1}^\infty\frac{\Lambda(n)}{\sqrt{n}}(\delta_{\log n}+\delta_{-\log n})+ 2\cosh(x/2)\,dx\]
is a tempered distribution. If this is true the measure $\mu$ is the Fourier transform of the tempered distribution 
\[\sum_{\gamma}\delta_{\gamma/2\pi}-2\vartheta'(2\pi t)\,dt,\]
where $\gamma$ run through the ordinates of the non-trivial zeros of $\zeta(s)$ repeated according to its multiplicity and $\vartheta(t)$ is the phase of the zeta function.
\end{theorem}
\begin{proof}
The measure $\mu$ is a Radon measure defined on $\R$. Therefore for $\varphi\in\mathcal{D}(\R)$ we have
\[\int\varphi(x)\mu(dx)=-\sum_{n=1}^\infty\frac{\Lambda(n)}{n^{1/2}}\bigl\{\varphi(\log n)+\varphi(-\log n)\bigr\}+2\int_{-\infty}^{+\infty}\cosh(x/2)\varphi(x)\,dx.\]

We have $\mu=d\alpha(x)$ where $\alpha(x)$ is a function of bounded variation in each interval and defined as $\alpha(x)=\mu(0,x]$ (for $x>0$). Hence
\[\alpha(x)=\sum_{\log n\le x}\frac{\Lambda(n)}{\sqrt{n}}-4\sinh(x/2).\]
For $x<0$ we have $\alpha(x)=-\alpha(-x)$ (for $x\not\in \Z$ and continuous to the right at every point). 

By Lemma \ref{E:201108-1} we obtain for $x>0$, and assuming the RH
\[\alpha(x)=\bigl(2e^{x/2}+\Orden(x^3)\bigr)-4\sinh(x/2)=\Orden(x^3).\]
The function $\alpha(x)$ is continuous in $\R$ except for jump discontinuities at the points $\pm\log(p^k)$. Hence, as distribution, it is the derivative of continuous function
that increases at most as $x^4$. Therefore $\mu$ is the second derivative of a continuous function that increases at most as $x^4$. By L. Schwartz \cite{S}*{Theorem VI in Chapter VII}, $\mu$ is a tempered distribution.

Assuming the Riemann hypothesis, we compute its Fourier transform.
The set of ordinates of the zeros of zeta is discrete, therefore the measure  $\sum_{\gamma}\delta_{\gamma/2\pi}$ is well defined as a Radon measure in $\R$. It is a tempered distribution according to \cite{S}*{Theorem VII, Ch. 7}, since the number of ceros $N(T)\ll T\log T$. Also $2\vartheta'(2\pi t)\,dt$ is tempered by the same reasoning. It follows that the measure $\nu:=\sum_{\gamma}\delta_{\gamma/2\pi}-2\vartheta'(2\pi t)\,dt$ is a tempered distribution.

There is a dense subspace $V$ of $\mathcal S(\R)$ of functions $f(t)$ satisfying the condition of Delsarte's Theorem \ref{T:Dels}. For example we may take as $V$ the linear combination of Hermite eigenfunctions of Fourier transform
\begin{equation}\label{E:eigenfunctions}
\Hermite_n(x)=\frac{2^{1/4}}{(2^n n!)^{1/2}}H_n(\sqrt{2\pi}\; x) e^{-\pi x^2}.
\end{equation}
Since we assume the Riemann hypothesis the non trivial zeros of zeta are 
$\rho=1/2+i\gamma$ and for  any function  $f\in V$, \eqref{E:Delsarte} may be written as  
\begin{equation}\label{E:Delsarte3}
\sum_\gamma f(\gamma)-F(0)-F(1)=-\frac{1}{2\pi}\sum_{n=1}^\infty \frac{\Lambda(n)}{\sqrt{n}}\bigl(\widehat f(\tfrac{\log n}{2\pi})+\widehat f (-\tfrac{\log n}{2\pi})\bigr)-\frac{1}{2\pi}\int_{-\infty}^\infty f(t)\Phi(\tfrac12+it)\,dt,
\end{equation}

Since $f\in V$
\[F(0)+F(1)=\int_{-\infty}^\infty \widehat f(-x)(e^{-\pi x}+e^{\pi x})\,dt=2\int_{-\infty}^\infty \widehat f(x)\cosh(\pi x)\,dx.\]
Since $\zeta(\frac12+it)=e^{-i\vartheta(t)}Z(t)$ we get 
\[\Phi(\tfrac12+it)=\frac{\zeta'(\frac12+it)}{\zeta(\frac12+it)}+\frac{\zeta'(\frac12-it)}{\zeta(\frac12-it)}=\frac{(-i\vartheta' Z+Z')e^{-i\vartheta}}{ie^{-i\vartheta}Z}-\frac{(i\vartheta' Z+Z')e^{i\vartheta}}{ie^{i\vartheta}Z}=-2\vartheta'.\]
Hence, using that $\vartheta'(t)$ is an even function
\begin{multline*}
\sum_\gamma f(\gamma)-\frac{1}{\pi}\int_{-\infty}^\infty f(t)\vartheta'(t)\,dt=\\-\frac{1}{2\pi}\sum_{n=1}^\infty \frac{\Lambda(n)}{\sqrt{n}}\bigl(\widehat f(\tfrac{\log n}{2\pi})+ \widehat f(-\tfrac{\log n}{2\pi})\bigr)+2\int_{-\infty}^\infty\widehat  f(x)\cosh(\pi x)\,dx.
\end{multline*}
With  $g(\frac{t}{2\pi})= f(t)$ we have  $\widehat f(x)=2\pi \widehat g(2\pi x)$ and therefore
\begin{multline*}
\sum_\gamma g(\tfrac{\gamma}{2\pi})-2\int_{-\infty}^\infty g(t)\vartheta'(2\pi t)\,dt=\\-\sum_{n=1}^\infty \frac{\Lambda(n)}{\sqrt{n}}\bigl(\widehat g(\log n)+ \widehat g(-\log n)\bigr)+2\int_{-\infty}^\infty\widehat  g(x)\cosh(x/2)\,dx.
\end{multline*}

That is 
\[\langle \nu,g\rangle=\langle \mu, \widehat g\rangle.\]
This is true for a dense set in $\mathcal{S}(\R)$, but since $\mu$ and $\nu$ are tempered distributions, the equality is also true for any $g\in\mathcal{S}(\R)$. 
By definition  $\nu=\widehat{\mu}$, since the measure are both even this is equivalent to $\mu=\widehat{\nu}$ i.e. 
\begin{equation}
\mathcal{F}\Bigl(\sum_\gamma\delta_{\gamma/2\pi}-2\vartheta'(2\pi t)\,dt\Bigr) =
-\sum_{n=1}^\infty\frac{\Lambda(n)}{\sqrt{n}}(\delta_{\log n}+\delta_{-\log n})+ 2\cosh(x/2)\,dx.
\end{equation}

To prove the last assertion, assume that the measure $\mu$ is a tempered distribution. For each $s$ with $\sigma>0$ the function $f_s(x)=\frac{1}{4\cosh(s x)}$ is in $\mathcal{S}(\R)$. It is also true that for any sequence $z_n\to0$ with $z_n\ne0$ we have 
\[\lim_{n\to\infty}\frac{f_{s+z_n}-f_s}{z_n}=f'_s, \qquad f'_s(x)=-\frac{ x\sinh(sx)}{4\cosh^2(s x)},\]
where the limit is in the usual topology of $\mathcal{S}(\R)$. 
It follows that the function $s\mapsto\langle f_s,\mu\rangle$ is a well defined holomorphic function for $\sigma>0$. 

For $\sigma> \frac12 $ we have 
\begin{align*}
\langle f_s,\mu\rangle&=-\sum_{n=1}^\infty \frac{\Lambda(n)}{\sqrt{n}}\frac{1}{2\cosh(s\log n)}+\int_{-\infty}^\infty\frac{\cosh(x/2)}{2\cosh(s x)}\,dx\\
&=-\sum_{n=2}^\infty \frac{\Lambda(n)}{n^{\frac12+s}}\sum_{k=0}^\infty(-1)^k n^{-2ks}+\frac{\pi}{2 s \cos(\frac{\pi}{4s})}\\
&=\sum_{k=0}^\infty(-1)^k\frac{\zeta'(\frac12+(2k+1)s)}{\zeta(\frac12+(2k+1)s)}+\frac{\pi}{2 s \cos(\frac{\pi}{4s})}
\end{align*}
The functions $(-1)^k\frac{\zeta'(\frac12+(2k+1)s)}{\zeta(\frac12+(2k+1)s)}$ and  $\frac{\pi}{2 s \cos(\frac{\pi}{4s})}$ both  have poles at  $s=\frac{1}{2(2k+1)}$ but its singular parts cancel out. 

Assuming $\mu$ is tempered implies that this function is holomorphic for $\sigma>0$. Any zero of $\zeta(s)$ to the right of $\sigma=\frac12$ will be a pole of this function situated at $\sigma>0$ and contradicting the hypothesis that $\mu$ is tempered.
\end{proof}

%This was written in November 2020 after my entry in the Blog of Imus
%https://institucional.us.es/blogimus/en/2020/11/frogs-birds-and-riemann-hypothesis/
%I put it now in arXiv since the paper by Felipe Gonçalves arXiv:2312.11185 writes the Guinand formula very near this one.

\end{document}